\documentstyle[amscd, leqno]{amsart}
\theoremstyle{plain}
\newtheorem{theorem}{Theorem}[section]
\newtheorem{proposition}[theorem]{Proposition}
\newtheorem{lemma}[theorem]{Lemma}

\theoremstyle{definition}

\theoremstyle{remark}
\newtheorem{remark}[theorem]{Remark}

\makeatletter

\@addtoreset{equation}{section}
\makeatother

\begin{document}
\title[]
{On the boundary behavior of the curvature of $L^{2}$-metrics}
\author{Ken-Ichi Yoshikawa}
\address{
Department of Mathematics,
Faculty of Science,
Kyoto University,
Kyoto 606-8502, JAPAN}
\email{yosikawa@@ms.u-tokyo.ac.jp}
\address{Korea Institute for Advanced Study,
Hoegiro 87, Dongdaemun-gu,
Seoul 130-722, KOREA}

\thanks{The author is partially supported by the Grants-in-Aid 
for Scientific Research (B) 19340016, JSPS}

\begin{abstract}
For one-parameter degenerations of compact K\"ahler manifolds,
we determine the asymptotic behavior of the first Chern form of  
the direct image of a Nakano semi-positive vector bundle twisted by the relative canonical bundle,
when the direct image is equipped with the $L^{2}$-metric.
\end{abstract}

\maketitle

\section
{Introduction}\label{sect:1}
\par
Let $X$ be a connected K\"ahler manifold of dimension $n+1$ with K\"ahler metric $h_{X}$
and let $S=\{s\in{\bf C};\,|s|<1\}$ be the unit disc. Set $S^{o}:=S\setminus\{0\}$.
Let $\pi\colon X\to S$ be a proper surjective holomorphic map with connected fibers. 
Let $\Sigma_{\pi}$ be the critical locus of $\pi$. We assume that $\pi(\Sigma_{\pi})=\{0\}$.  
We set $X_{s}=\pi^{-1}(s)$ for $s\in S$. Then $X_{s}$ is non-singular for $s\in S^{o}$.
Let $\omega_{X}=\Omega_{X}^{n+1}$ be the canonical bundle of $X$ 
and let $\omega_{X/S}=\Omega_{X}^{n+1}\otimes(\pi^{*}\Omega^{1}_{S})^{-1}$
be the relative canonical bundle of $\pi\colon X\to S$. 
The K\"ahler metric $h_{X}$ induces a Hermitian metric
$h_{X/S}$ on $TX/S=\ker\pi_{*}|_{X\setminus\Sigma_{\pi}}$, and $h_{X/S}$ induces a Hermitian metric
$h_{\omega_{X/S}}$ on $\omega_{X/S}$.
\par
Let $\xi\to X$ be a holomorphic vector bundle on $X$ equipped with a Hermitian metric $h_{\xi}$.
We write $\omega_{X/S}(\xi)=\omega_{X/S}\otimes\xi$.
In this note, we assume that $(\xi,h_{\xi})$ is a Nakano semi-positive vector bundle on $X$.
Namely, if $R^{\xi}$ denotes the curvature form of $(\xi,h_{\xi})$ with respect to the holomorphic Hermitian
connection, then the Hermitian form $h_{\xi}(\sqrt{-1}R^{\xi}(\cdot),\cdot)$ on the holomorphic vector bundle 
$TX\otimes\xi$ is semi-positive.
Since $\dim S=1$ and since $(\xi,h_{\xi})$ is Nakano semi-positive, 
all direct image sheaves $R^{q}\pi_{*}\omega_{X/S}(\xi)$ are locally free by \cite{Takegoshi95}. 
By the fiberwise Hodge theory, $R^{q}\pi_{*}\omega_{X/S}(\xi)$ is equipped with
the $L^{2}$-metric $h_{L^{2}}$ with respect to $h_{X/S}$ and $h_{\omega_{X/S}}\otimes h_{\xi}$.
By Berndtsson \cite{Berndtsson09} and Mourougane-Takayama \cite{MourouganeTakayama08}, 
the holomorphic Hermitian vector bundle $(R^{q}\pi_{*}\omega_{X/S}(\xi),h_{L^{2}})$ 
is again Nakano semi-positive on $S^{o}$.
By Mourougane-Takayama \cite{MourouganeTakayama09}, $h_{L^{2}}$ induces a singular Hermitian metric
with semi-positive curvature current on the tautological quotient bundle over the projective-space bundle 
${\bf P}(R^{q}\pi_{*}\omega_{X/S}(\xi))$. (We remark that there is no restrictions of the dimension of the base
space $S$ in the works \cite{Berndtsson09}, \cite{MourouganeTakayama08}, \cite{MourouganeTakayama09}.)
\par
After these results, one of the natural problems to be considered 
is the quantitative estimates for the singularities of the $L^{2}$-metric and its curvature.
In \cite{Yoshikawa10}, we gave a formula for the singularity of the $L^{2}$-metric on 
$R^{q}\pi_{*}\omega_{X/S}(\xi)$ (cf. Sect.\,\ref{sect:2}).
As a consequence, if $\sigma_{q}$ is a nowhere vanishing holomorphic section of
$\det R^{q}\pi_{*}\omega_{X/S}(\xi)$, then there exist a rational number $a_{q}\in{\bf Q}$,
an integer $\ell_{q}\geq0$ and a real number $c_{q}$ such that (cf. \cite[Th.\,6.8]{Yoshikawa10})
$$
\log\|\sigma_{q}(s)\|_{L^{2}}^{2}
=
a_{q}\,\log|s|^{2}+\ell_{q}\log(-\log|s|^{2})+c_{q}+O(1/\log|s|)
\qquad
(s\to0).
$$
\par
In this note, we study the boundary behavior of the curvature of the holomorphic Hermitian vector bundle
$(R^{q}\pi_{*}\omega_{X/S}(\xi),h_{L^{2}})$ as an application of the description of the singularity 
of the $L^{2}$-metric $h_{L^{2}}$ given in \cite{Yoshikawa10}. 
In this sense, this note is a supplement to the article \cite{Yoshikawa10}.
\par
Let us state our results.
Let ${\mathcal R}(s)\,ds\wedge d\bar{s}$ be the curvature form of $R^{q}\pi_{*}\omega_{X/S}(\xi)$
with respect to the holomorphic Hermitian connection associated to $h_{L^{2}}$.
By the Nakano semi-positivity \cite{Berndtsson09}, \cite{MourouganeTakayama08},
$\sqrt{-1}{\mathcal R}(s)$ is a semi-positive Hermitian endomorphism on the Hermitian bundle 
$(R^{q}\pi_{*}\omega_{X/S}(\xi),h_{L^{2}})$ on $S^{o}$.

\begin{theorem}\label{MainTheorem}
The curvature form ${\mathcal R}(s)\,ds\wedge d\bar{s}$ has Poincar\'e growth near $0\in S$. 
Namely, there exists a constant $C>0$ such that the following inequality of Hermitian endomorphisms holds
for all $s\in S^{o}$
$$
0\leq
\sqrt{-1}{\mathcal R}(s)
\leq
\frac{C}{|s|^{2}(\log|s|)^{2}}\,{\rm Id}_{R^{q}\pi_{*}\omega_{X/S}(\xi)}.
$$
Moreover, the Chern form $c_{1}(R^{q}\pi_{*}\omega_{X/S}(\xi),h_{L^{2}})$ 
has the following asymptotic behavior as $s\to0$:
$$
c_{1}(R^{q}\pi_{*}\omega_{X/S}(\xi),h_{L^{2}})
=
\left\{
\frac{\ell_{q}}{|s|^{2}(\log|s|)^{2}}
+
O\left(\frac{1}{|s|^{2}(\log|s|)^{3}}\right)
\right\}
\sqrt{-1}\,ds\wedge d\bar{s}.
$$
\end{theorem}

Considering the trivial example $X=M\times S$, $\xi={\mathcal O}_{X}$, $\pi={\rm pr}_{2}$, where $M$
is a compact K\"ahler manifold, we can not expect any lower bound of $\sqrt{-1}{\mathcal R}(s)$ (resp.
$c_{1}(R^{q}\pi_{*}\omega_{X/S}(\xi),h_{L^{2}})$) by a non-zero semi-positive Hermitian endomorphism
(resp. real $(1,1)$-form). 
We remark that, when $X_{0}$ is reduced and has only canonical singularities, 
then we get a better estimate (cf. Sect.\,\ref{sect:5}).
\par
As an application of Theorem~\ref{MainTheorem}, we get an estimate for the complex Hessian of analytic torsion.
Set $X_{s}:=\pi^{-1}(s)$ and $\xi_{s}:=\xi|_{X_{s}}$ for $s\in S$.
Let $\omega_{X_{s}}$ be the canonical line bundle of $X_{s}$ 
and let $h_{\omega_{X_{s}}}$ be the Hermitian metric on $\omega_{X_{s}}$ induced from $h_{X}$.
For $s\in S^{o}$, let $\tau(X_{s},\omega_{X_{s}}(\xi_{s}))$ be the analytic torsion \cite{RaySinger73},
\cite{BGS88} of the holomorphic Hermitian vector bundle
$(\xi_{s}\otimes\omega_{X_{s}},h_{\xi}|_{X_{s}}\otimes h_{\omega_{X_{s}}})$ 
on the compact K\"ahler manifold $(X_{s},h_{X}|_{X_{s}})$. 
Let $\log\tau(X/S,\omega_{X/S}(\xi))$ be the function defined as
$$
\log\tau(X/S,\omega_{X/S}(\xi))(s):=\log\tau(X_{s},\omega_{X_{s}}(\xi_{s})),
\qquad
s\in S^{o}.
$$
By Bismut-Gillet-Soul\'e \cite{BGS88}, $\log\tau(X/S,\omega_{X/S}(\xi))$ is a $C^{\infty}$ function on $S^{o}$. 
Moreover, under certain algebraicity assumption of the family $\pi\colon X\to S$ and the vector bundle $\xi$,
there exist by \cite{Yoshikawa10} constants 
$\alpha\in{\bf Q}$, $\beta\in{\bf Z}$, $\gamma\in{\bf R}$ such that 
$$
\log\tau(X/S,\omega_{X/S}(\xi))(s)
=
\alpha\,\log|s|^{2}-(\sum_{q\geq0}(-1)^{q}\ell_{q})\,\log(-\log|s|^{2})+\gamma+O(1/\log|s|)
$$
as $s\to0$.
By this asymptotic expansion, it is reasonable to expect that the complex Hessian of analytic torsion
has a similar behavior to the Poincar\'e metric on $S^{o}$.

\begin{theorem}\label{Theorem:HessianTorsion}
The complex Hessian $\partial_{s\bar{s}}\log\tau(X/S,\omega_{X/S}(\xi))$ has the following 
asymptotic behavior as $s\to0$:
$$
\partial_{s\bar{s}}\log\tau(X/S,\omega_{X/S}(\xi))
=
\frac{\sum_{q\geq0}(-1)^{q}\ell_{q}}{|s|^{2}(\log|s|)^{2}}
+
O\left(\frac{1}{|s|^{2}(\log|s|)^{3}}\right).
$$
\end{theorem}

This note is organized as follows.
In Sect.\,\ref{sect:2}, we recall the structure of the singularity of the $L^{2}$-metric $h_{L^{2}}$
on $R^{q}\pi_{*}\omega_{X/S}(\xi)$. 
In Sect.\,\ref{sect:3}, we prove some technical lemmas used in the proof of Theorem~\ref{MainTheorem}.
In Sect.\,\ref{sect:4}, we prove Theorems~\ref{MainTheorem} and \ref{Theorem:HessianTorsion}.
In Sect.\,\ref{sect:5}, we study the case where $X_{0}$ has only canonical singularities. 
\par
Throughout this note, we keep the notation and the assumptions in Sect.\,\ref{sect:1}.

\section
{The singularity of the $L^{2}$-metric}\label{sect:2}
\par

\subsection
{The structure of the singularity of the $L^{2}$-metric}\label{subsect:2.1}
\par
Let $\kappa_{\mathcal X}$ be the K\"ahler form of $h_{\mathcal X}$. 
In the rest of this note, we assume that 
{\em $(\xi,h_{\xi})$ is Nakano semi-positive on $X$ and that $(S,0)\cong(\varDelta,0)$.}
By  \cite[Th.\,6.5 (i)]{Takegoshi95}, $R^{q}\pi_{*}\omega_{X/S}(\xi)$ is locally free on $S$. 
By shrinking $S$ if necessary, we may also assume that $R^{q}\pi_{*}\omega_{X/S}(\xi)$ is a free
${\mathcal O}_{S}$-module on $S$.
Let $\rho_{q}\in{\bf Z}_{\geq0}$ be the rank of $R^{q}\pi_{*}\Omega_{X}^{n+1}(\xi)$ 
as a free ${\mathcal O}_{S}$-module on $S$. 
Let $\{\psi_{1},\ldots,\psi_{\rho_{q}}\}\subset H^{0}(S,R^{q}\pi_{*}\omega_{X/S}(\xi))$
be a free basis of the locally free sheaf $R^{q}\pi_{*}\omega_{X/S}(\xi)$ on $S$.
\par
Let $T$ be another unit disc.
By the semistable reduction theorem \cite[Chap II]{Mumford77}, there exists a positive integer $\nu>0$
such that the family $X\times_{S}T\to T$ induced from $\pi\colon X\to S$
by the map $\mu\colon T\to S$, $\mu(t)=t^{\nu}$, admits a semistable model. 
Namely, there is a resolution $r\colon Y\to X\times_{S}T$ such that the family 
$f:={\rm pr}_{2}\circ r\colon Y\to T$ is semistable, i.e., 
$Y_{0}:=f^{-1}(0)$ is a {\em reduced} normal crossing divisor of $Y$. We fix such an integer $\nu>0$.
Let ${\rm Herm}(r)$ be the set of $r\times r$-Hermitian matrices.

\begin{theorem}\label{Theorem2.1}
By choosing a basis $\{\psi_{1},\ldots,\psi_{\rho_{q}}\}$ of $R^{q}\pi_{*}\omega_{X/S}(\xi)$ 
as a free ${\mathcal O}_{S}$-module appropriately, the $\rho_{q}\times\rho_{q}$-Hermitian matrix 
$$
G(s):=(h_{L^{2}}(\psi_{\alpha}|_{X_{s}},\psi_{\beta}|_{X_{s}}))
$$
has the following expression
$$
G(t^{\nu})=D(t)\cdot H(t)\cdot\overline{D(t)},
\qquad
D(t)={\rm diag}(t^{-e_{1}},\ldots,t^{-e_{\rho_{q}}}).
$$
Here $e_{1},\ldots,e_{\rho_{q}}\geq0$ are integers and the Hermitian matrix $H(t)$ has the following structure:
There exist $A_{m}(t)\in C^{\infty}(T,{\rm Herm}(\rho_{q}))$, $0\leq m\leq n$, with
$$
H(t)=\sum_{m=0}^{n}A_{m}(t)\,(\log |t|^{2})^{m}.
$$
Moreover, by defining the real-valued functions $a_{m}(t)\in C^{\infty}(T)$, $0\leq m\leq n\rho_{q}$ as
$$
\det H(t)=\sum_{m=0}^{n\rho_{q}}a_{m}(t)\,(\log|t|^{2})^{m},
$$
one has $a_{m}(0)\not=0$ for some $0\leq m\leq n\rho_{q}$.
\end{theorem}

\begin{pf}
See \cite[Th.\,6.8 and Lemmas 6.3 and 6.4]{Yoshikawa10}.
\end{pf}

\begin{remark}\label{Remark2.2}
The meaning of the Hermitian matrix $H(t)$ and the diagonal matrix $D(t)$  is explained as follows. 
Let $F\colon Y\to X$ be the map defined as the composition of $r\colon Y\to X\times_{S}T$ and
${\rm pr}_{1}\colon X\times_{S}T\to X$. Then $h_{Y}:=r^{*}(h_{X}+dt\otimes d\bar{t})$ is a K\"ahler metric
on $Y\setminus Y_{0}$. There is a basis $\{\theta_{1},\ldots,\theta_{\rho_{q}}\}$ of
$R^{q}f_{*}\omega_{Y/Y}(F^{*}\xi)$ such that 
$H(t)=(H_{\alpha\bar{\beta}}(t))$,
$H_{\alpha\bar{\beta}}(t)=(\mu^{*}h_{L^{2}})(\theta_{\alpha},\theta_{\beta})$, 
where $\mu^{*}h_{L^{2}}$ is the $L^{2}$-metric
on $R^{q}f_{*}\omega_{Y/T}(F^{*}\xi)$ with respect to $h_{Y}$ and $F^{*}h_{\xi}$.
By \cite[Lemma 3.3]{MourouganeTakayama09}, $R^{q}f_{*}\omega_{Y/T}(F^{*}\xi)$ is regarded as
a subsheaf of $\mu^{*}R^{q}f_{*}\omega_{X/S}(\xi)$. Then the relation between the two basis
$\{\theta_{1},\ldots,\theta_{\rho_{q}}\}$ and $\{\mu^{*}\psi_{1},\ldots,\mu^{*}\psi_{\rho_{q}}\}$
is given by $D(t)$, i.e., $\theta_{\alpha}=t^{e_{\alpha}}\mu^{*}\psi_{\alpha}$.
Moreover, by \cite[Lemma 4.2]{MourouganeTakayama09}, 
$\mu^{*}h_{L^{2}}|_{T^{o}}$ is indeed the pull-back of the $L^{2}$-metric $h_{L^{2}}|_{S^{o}}$ via $\mu$, 
where $T^{o}:=T\setminus\{0\}$, which implies the relation
$G(\mu(t))=D(t)H(t)\overline{D(t)}$.
\par
The proof of Theorem~\ref{Theorem2.1} heavily relies on a theorem of Barlet \cite[Th.\,4 bis.]{Barlet82}.
This is the major reason why we need the assumption $\dim S=1$.
\end{remark}

\subsection
{A Hodge theoretic proof of Theorem~\ref{Theorem2.1} for a trivial line bundle}\label{subsect:2.2}
\par
Assume that $(\xi,h_{\xi})$ is a trivial Hermitian line bundle on $X$, 
that $\pi\colon X\to S$ is a family of polarized projective manifolds with unipotent monodromy 
and that the K\"ahler class of $h_{X}$ is the first Chern class of an ample line bundle on $X$.
We see that the expansion in Theorem~\ref{Theorem2.1} follows from the nilpotent orbit theorem of 
Schmid \cite{Schmid73} in this case. 
\par
Let $\kappa_{X}$ be the K\"ahler class of $h_{X}$. By assumption,
there is a very ample line bundle $L$ on $X$ with $[\kappa_{X}]=c_{1}(L)/N$. 
Replacing $\kappa_{X}$ by $N\kappa_{X}$ if necessary, we may assume that $L$ is very ample.
Let $H_{1},\ldots,H_{n}\in |L|$ be sufficiently generic hyperplane sections such that
the following hold for all $0\leq k\leq n$ after shrinking $S$ if necessary:
\begin{itemize}
\item[(i)]
$X\cap H_{1}\cap\cdots\cap H_{k}$ is a complex manifold of dimension $n-k+1$.
\item[(ii)]
The restriction of $\pi$ to $X\cap H_{1}\cap\cdots\cap H_{k}$ is a flat holomorphic map from 
$X\cap H_{1}\cap\cdots\cap H_{k}$ to $S$. 
\item[(iii)]
$X_{s}\cap H_{1}\cap\cdots\cap H_{k}$ is a projective manifold of dimension $n-k$ for $s\in S^{o}$.
\end{itemize}
We set $X_{s}^{(k)}:=(\pi^{(k)})^{-1}(s)=X_{s}\cap H_{1}\cap\cdots\cap H_{k}$ for $s\in S$. 
\par
Let $\{\psi_{1},\ldots,\psi_{\rho_{q}}\}\subset H^{0}(S,R^{q}\pi_{*}\omega_{X/S})$
be a free basis of the locally free sheaf $R^{q}\pi_{*}\omega_{X/S}$ on $S$. 
There exists 
$\Psi_{1},\ldots,\Psi_{\rho_{q}}\in H^{0}(X,\Omega_{X}^{n+1-q})$ by \cite[Th.\,5.2]{Takegoshi95} 
(after schrinking $S$ if necessary) such that
$$
\psi_{\alpha}=[(\Psi_{\alpha}\wedge\kappa_{X}^{q})\otimes(\pi^{*}ds)^{-1}],
\qquad
\pi^{*}ds\wedge\Psi_{\alpha}=0.
$$
By the condition $\pi^{*}ds\wedge\Psi_{\alpha}=0$, there exist relative holomorphic differentials 
$\psi'_{\alpha}\in H^{0}(X\setminus\Sigma_{\pi},\Omega_{X/S}^{n-q})$ such that
$\Psi_{\alpha}=\psi'_{\alpha}\wedge\pi^{*}ds$. 
Then the harmonic representative of the cohomology class $\psi|_{X_{s}}$ is given by
$\psi'_{\alpha}\wedge\kappa_{X}|_{X_{s}}$.
Since $\kappa_{X}=c_{1}(L)$, we get
$$
h_{L^{2}}(\psi_{\alpha},\psi_{\beta})(s)
=
i^{(n-q)^{2}}\int_{X_{s}}\psi'_{\alpha}\wedge\overline{\psi'}_{\beta}\wedge\kappa_{X}^{q}|_{X_{s}}
=
i^{(n-q)^{2}}\int_{X_{s}^{(q)}}\psi'_{\alpha}\wedge\overline{\psi'}_{\beta}|_{X^{(q)}_{s}}.
$$
Hence Theorem~\ref{Theorem2.1} is reduced to the case $q=0$.
In the case $q=0$, Theorem~\ref{Theorem2.1} is a consequence of Fujita's estimate \cite[1.12]{Fujita78}
and the following:

\begin{lemma}\label{Lemma2.3}
For $\varphi,\psi\in H^{0}(X,\Omega_{X}^{n+1})$,
there exist $a_{m}(s)\in C^{\omega}(S)$, $0\leq m\leq n$ such that
$$
\pi_{*}(\varphi\wedge\overline{\psi})(s)
=
\sum_{m=0}^{n}(\log|s|^{2})^{m}a_{m}(s)\,ds\wedge d\bar{s}.
$$
In particular,
$h_{L^{2}}(\varphi\otimes(\pi^{*}ds)^{-1}|_{X_{s}},\psi\otimes(\pi^{*}ds)^{-1}|_{X_{s}})=
i^{n^{2}}\sum_{m=0}^{n}(\log|s|^{2})^{m}a_{m}(s)$.
\end{lemma}

\begin{pf}
Fix ${\frak o}\in S^{o}$. Let $\gamma\in GL(H^{n}(X_{\frak o},{\bf C}))$ be the monodromy.
By assumption, $\gamma$ is unipotent.
Set ${\bf H}:=R^{n}\pi_{*}{\bf C}\otimes_{\bf C}{\mathcal O}_{S^{o}}$, which is equipped with 
the Gauss-Manin connection. Let $\{v_{1},\ldots,v_{m}\}$ be a basis of $H^{n}(X_{\frak o},{\bf C})$.
Since $\gamma$ is unipotent, there exists a nilpotent $N\in{\rm End}(H^{n}(X_{\frak o},{\bf C}))$ 
such that $\gamma=\exp(N)$. 
Let $p\colon\widetilde{S^{o}}\ni z\to\exp(2\pi iz)\in S^{o}$ be the universal covering.
Since ${\bf H}$ is flat, $v_{k}$ extend to flat sections 
${\bf v}_{k}\in\Gamma(\widetilde{S^{o}},p^{*}{\bf H})$, which induces an isomorphism
$p^{*}{\bf H}\cong{\mathcal O}_{\widetilde{S^{o}}}\otimes_{\bf C}H^{n}(X_{\frak o},{\bf C})$.
Under this trivialization, we have ${\bf v}_{k}(z+1)=\gamma\cdot{\bf v}_{k}(z)$.
We define ${\bf s}_{k}(\exp 2\pi iz):=\exp\left(-z\,N\right)\,{\bf v}_{k}(z)$.
Then ${\bf s}_{1},\ldots,{\bf s}_{m}\in\Gamma(\widetilde{S^{o}},p^{*}{\bf H})$ descend to 
single-valued holomorphic frame fields of ${\bf H}$. 
The canonical extension of ${\bf H}$ is the locally free sheaf on $S$ defined as 
$\overline{\bf H}:={\mathcal O}_{S}\,{\bf s}_{1}\oplus\cdots\oplus{\mathcal O}_{S}\,{\bf s}_{m}$.
Set ${\bf F}^{n}:=\pi_{*}\Omega^{n}_{X/S}|_{S^{o}}\subset{\bf H}$. By \cite[p.\,235]{Schmid73},
${\bf F}^{n}$ extends to a subbundle $\overline{\bf F}^{n}\subset\overline{\bf H}$. 
\par
There exists 
$\varphi',\psi'\in H^{0}(X\setminus X_{0},\Omega_{X/S}^{n}|_{X\setminus X_{0}})$ such that
$\varphi=\pi^{*}ds\wedge\varphi'$ and $\psi=\pi^{*}ds\wedge\psi'$ on $X\setminus X_{0}$. 
Then $\varphi'$ and $\psi'$ are identified with 
$\varphi\otimes(\pi^{*}ds)^{-1},\psi\otimes(\pi^{*}ds)^{-1}\in H^{0}(X,\omega_{X/S})$, respectively. 
Since ${\bf F}^{n}\subset{\bf H}$, there exist $b_{k}(t),c_{k}(t)\in{\mathcal O}(S^{o})$ such that
$[\varphi'|_{X_{s}}]=\sum_{k=1}^{m}b_{k}(s)\,{\bf s}_{k}(s)$ and
$[\psi'|_{X_{s}}]=\sum_{k=1}^{m}c_{k}(s)\,{\bf s}_{k}(s)$. 
Since $\pi_{*}\omega_{X/S}=\overline{F}^{n}$ by Kawamata \cite[Lemma 1]{Kawamata82}, 
we get $b_{k}(s),c_{k}(s)\in{\mathcal O}(S)$. Then
$$
\pi_{*}(\varphi\wedge\overline{\psi})(s)
=
\{\int_{X_{s}}\varphi'\wedge\overline{\psi'}\}\,ds\wedge d\bar{s}
=
\{\int_{X_{s}}\sum_{j=1}^{m}b_{j}(s){\bf s}_{j}(s)\wedge
\sum_{k=1}^{m}\overline{c_{k}(s){\bf s}_{k}(s)}\}\,ds\wedge d\bar{s}.
$$
Substituting 
${\bf s}_{k}(s)=\exp\left(-z\,N\right){\bf v}_{k}(z)=\sum_{0\leq m\leq n}\frac{(-z)^{m}}{m!}N^{m}{\bf v}_{k}(z)$, 
we get
$$
\pi_{*}(\varphi\wedge\overline{\psi})(s)
=
\{\sum_{j,k=1}^{m}b_{j}(s)\overline{c_{k}(s)}
\sum_{0\leq a,b\leq n}\frac{(-1)^{a+b}}{a!b!}z^{a}\bar{z}^{b}
C_{a,b}^{j,k}\}\,
ds\wedge d\bar{s},
$$
where $z=\frac{1}{2\pi i}\log s$ and
$C_{a,b}^{j,k}=\int_{X_{\frak o}}(N^{a}v_{j})\wedge(N^{b}\overline{v_{k}})$.
Since $\pi_{*}(\varphi\wedge\overline{\psi})$ is single-valued, so is the expression 
$\sum_{a+b=m}\frac{(-1)^{a+b}}{a!b!}z^{a}\bar{z}^{b}C_{a,b}^{j,k}$. 
As a result, there exists a constant $C^{j,k}_{m}\in{\bf C}$ such that
$\sum_{a+b=m}\frac{(-1)^{a+b}}{a!b!}z^{a}\bar{z}^{b}C_{a,b}^{j,k}=C^{j,k}_{m}(\log|s|^{2})^{m}$.
Setting $a_{m}(s):=\sum_{j,k=1}^{m}C^{j,k}_{m}b_{j}(s)\overline{c_{k}(s)}$, we get the result.
\end{pf}

\begin{remark}\label{Remark2.4}
In the proof of Theorem~\ref{Theorem2.1}, 
the role of the nilpotent orbit theorem is played by Barlet's theorem \cite[Th.\,4bis.]{Barlet82}
on the asymptotic expansion of fiber integrals associated to the function 
$f(z)=z_{0}\cdots z_{n}$ near the origin. See \cite[Sects.\,6.3 and 6.4]{Yoshikawa10} for more details.
\end{remark}

\section
{Some technical lemmas}\label{sect:3}
\par
We denote by $C^{\infty}_{\bf R}(T)$ the set of {\em real-valued} $C^{\infty}$ functions on $T$.

\begin{lemma}\label{Lemma3.1}
Let $\varphi(t)\in C^{\infty}_{\bf R}(T)$ and let $r\in{\bf Q}$ and $\ell\in{\bf Z}$.
Set $h(t):=|t|^{2r}(\log|t|^{2})^{\ell}\varphi(t)$. Then the following identities hold:
$$
\begin{array}{lll}
(1)&
\partial_{t}h(t)
&=
\left(
\frac{r}{t}+\frac{\ell}{t(\log|t|^{2})}+\frac{\partial_{t}\varphi(t)}{\varphi(t)}
\right)
h(t),
\\
(2)&
\partial_{t\bar{t}}h(t)
&=
\left(
-\frac{\ell}{|t|^{2}(\log|t|^{2})}+\frac{\partial_{t\bar{t}}\varphi(t)}{\varphi(t)}
-\frac{|\partial_{t}\varphi(t)|^{2}}{\varphi(t)^{2}}
+
\left|
\frac{r}{t}+\frac{\ell}{t(\log|t|^{2})}+\frac{\partial_{t}\varphi(t)}{\varphi(t)}
\right|^{2}
\right)
h(t).
\end{array}
$$
\end{lemma}

\begin{pf}
The proof is elementary and is left to the reader.
\end{pf}

\begin{lemma}\label{Lemma3.2}
Let $I$ be a finite set. For $i\in I$, let $r_{i}\in{\bf Q}$, $\ell_{i}\in{\bf Z}$ and $\varphi_{i}(t)\in C^{\infty}_{\bf R}(T)$. 
Set $g_{i}(t):=|t|^{2r_{i}}(\log|t|^{2})^{\ell_{i}}\varphi_{i}(t)$ for $i\in I$ and $g(t):=\sum_{i\in I}g_{i}(t)$. 
\begin{itemize}
\item[(1)]
If $g(t)>0$ on $T^{o}$, then the following equalities of functions on $T^{o}$ hold:
$$
\partial_{t}\log g
=
\sum_{i\in I}
\left(
\frac{r_{i}}{t}+\frac{\ell_{i}}{t(\log|t|^{2})}+\frac{\partial_{t}\varphi_{i}}{\varphi_{i}}
\right)
\frac{g_{i}}{g},
$$
$$
\begin{aligned}
\partial_{t\bar{t}}\log g
&=
-\frac{1}{2}\sum_{i,j}\frac{\ell_{i}+\ell_{j}}{|t|^{2}(\log|t|^{2})^{2}}
\cdot
\frac{g_{i}g_{j}}{g^{2}}
+
\frac{1}{2}\sum_{i,j}
\left|
\frac{r_{i}-r_{j}}{t}+\frac{\ell_{i}-\ell_{j}}{t(\log|t|^{2})}\right|^{2}
\cdot
\frac{g_{i}g_{j}}{g^{2}}
\\
&\quad
+
\frac{1}{2}\sum_{i,j}
\left(
\frac{\partial_{t\bar{t}}\varphi_{i}}{\varphi_{i}}+\frac{\partial_{t\bar{t}}\varphi_{j}}{\varphi_{j}}
\right)
\frac{g_{i}g_{j}}{g^{2}}
+
\frac{1}{2}\sum_{i,j}
\Re\left(
\frac{\partial_{t}\varphi_{i}}{\varphi_{i}}\cdot\overline{\frac{\partial_{t}\varphi_{j}}{\varphi_{j}}}
\right)
\frac{g_{i}g_{j}}{g^{2}}
\\
&\quad
+
\sum_{i<j}
\Re\left\{
\left(\frac{r_{i}-r_{j}}{t}+\frac{\ell_{i}-\ell_{j}}{t(\log|t|^{2})}\right)
\cdot
\left(\overline{\frac{\partial_{t}\varphi_{i}}{\varphi_{i}}-\frac{\partial_{t}\varphi_{j}}{\varphi_{j}}}\right)
\right\}
\frac{g_{i}g_{j}}{g^{2}}.
\end{aligned}
$$
\item[(2)]
If $r_{i}\geq0$ and $0\leq \ell_{i}\leq N$ for all $i\in I$, then as $t\to0$
$$
\begin{aligned}
\partial_{t\bar{t}}\log g
&=
-\frac{1}{2}\sum_{i,j}\frac{\ell_{i}+\ell_{j}}{|t|^{2}(\log|t|^{2})^{2}}
\cdot
\frac{g_{i}g_{j}}{g^{2}}
+
\frac{1}{2}\sum_{i,j}
\left|
\frac{r_{i}-r_{j}}{t}+\frac{\ell_{i}-\ell_{j}}{t(\log|t|^{2})}\right|^{2}
\cdot
\frac{g_{i}g_{j}}{g^{2}}
\\
&\qquad
+
O\left(\frac{(-\log|t|)^{2N}}{|t|g(t)^{2}}
\right).
\end{aligned}
$$
\end{itemize}
\end{lemma}

\begin{pf}
{\bf (1) }
The first equality of (1) follows from Lemma~\ref{Lemma3.1} (1).
Since
$$
\partial_{t\bar{t}}g(t)
=
\sum_{i\in I}
\left(
-\frac{\ell_{i}}{|t|^{2}(\log|t|^{2})}
+
\frac{\partial_{t\bar{t}}\varphi_{i}}{\varphi_{i}}
-\frac{|\partial_{t}\varphi_{i}|^{2}}{\varphi_{i}^{2}}
+
\left|
\frac{r_{i}}{t}+\frac{\ell_{i}}{t(\log|t|^{2})}+\frac{\partial_{t}\varphi_{i}}{\varphi_{i}}
\right|^{2}
\right)
g_{i}
$$
by Lemma~\ref{Lemma3.1} (2), we get
\begin{equation}\label{eqn:(4.1)}
\begin{aligned}
g\partial_{t\bar{t}}g
&=
\sum_{i,j\in I}g_{j}
\left(
-\frac{\ell_{i}}{|t|^{2}(\log|t|^{2})}
+
\frac{\partial_{t\bar{t}}\varphi_{i}}{\varphi_{i}}
-
\frac{|\partial_{t}\varphi_{i}|^{2}}{\varphi_{i}^{2}}
+
\left|
\frac{r_{i}}{t}+\frac{\ell_{i}}{t(\log|t|^{2})}+\frac{\partial_{t}\varphi_{i}}{\varphi_{i}}
\right|^{2}
\right)
g_{i}
\\
&=
\frac{1}{2}\sum_{i,j\in I}
\left(
-\frac{\ell_{i}+\ell_{j}}{|t|^{2}(\log|t|^{2})}
+
\frac{\partial_{t\bar{t}}\varphi_{i}}{\varphi_{i}}
+
\frac{\partial_{t\bar{t}}\varphi_{j}}{\varphi_{j}}
-
\frac{|\partial_{t}\varphi_{i}|^{2}}{\varphi_{i}^{2}}
-
\frac{|\partial_{t}\varphi_{j}|^{2}}{\varphi_{j}^{2}}
\right.
\\
&\qquad\qquad\quad
\left.
+
\left|
\frac{r_{i}}{t}+\frac{\ell_{i}}{t(\log|t|^{2})}+\frac{\partial_{t}\varphi_{i}}{\varphi_{i}}
\right|^{2}
+
\left|
\frac{r_{j}}{t}+\frac{\ell_{j}}{t(\log|t|^{2})}+\frac{\partial_{t}\varphi_{j}}{\varphi_{j}}
\right|^{2}
\right)
g_{i}g_{j}.
\end{aligned}
\end{equation}
By the first equality of Lemma~\ref{Lemma3.2} (1), we get
\begin{equation}\label{eqn:(4.2)}
\begin{aligned}
|\partial_{t}g|^{2}
&=
\sum_{i,j\in I}
\left(
\frac{r_{i}}{t}+\frac{\ell_{i}}{t(\log|t|^{2})}+\frac{\partial_{t}\varphi_{i}}{\varphi_{i}}
\right)
\left(
\overline{\frac{r_{j}}{t}+\frac{\ell_{j}}{t(\log|t|^{2})}+\frac{\partial_{t}\varphi_{j}}{\varphi_{j}}}
\right)
g_{i}g_{j}
\\
&=
\sum_{i,j\in I}
\Re\left(
\frac{r_{i}}{t}+\frac{\ell_{i}}{t(\log|t|^{2})}+\frac{\partial_{t}\varphi_{i}}{\varphi_{i}}
\right)
\left(
\overline{\frac{r_{j}}{t}+\frac{\ell_{j}}{t(\log|t|^{2})}+\frac{\partial_{t}\varphi_{j}}{\varphi_{j}}}
\right)
g_{i}g_{j}.
\end{aligned}
\end{equation}
By \eqref{eqn:(4.1)} and \eqref{eqn:(4.2)}, we get
\begin{equation}\label{eqn:(4.3)}
\begin{aligned}
g\partial_{t\bar{t}}g-|\partial_{t}g|^{2}
&=
\frac{1}{2}\sum_{i,j\in I}
\left(
-\frac{\ell_{i}+\ell_{j}}{|t|^{2}(\log|t|^{2})}
+
\frac{\partial_{t\bar{t}}\varphi_{i}}{\varphi_{i}}
+
\frac{\partial_{t\bar{t}}\varphi_{j}}{\varphi_{j}}
-
\frac{|\partial_{t}\varphi_{i}|^{2}}{\varphi_{i}^{2}}
-
\frac{|\partial_{t}\varphi_{j}|^{2}}{\varphi_{j}^{2}}
\right.
\\
&\quad
\left.
+
\left|
\left(
\frac{r_{i}}{t}+\frac{\ell_{i}}{t(\log|t|^{2})}+\frac{\partial_{t}\varphi_{i}}{\varphi_{i}}
\right)
-
\left(
\frac{r_{j}}{t}+\frac{\ell_{j}}{t(\log|t|^{2})}+\frac{\partial_{t}\varphi_{j}}{\varphi_{j}}
\right)
\right|^{2}
\right)
g_{i}g_{j}.
\end{aligned}
\end{equation}
Since $\partial_{t\bar{t}}\log g=(g\partial_{t\bar{t}}g-|\partial_{t}g|^{2})/g^{2}$,
the second equality of Lemma~\ref{Lemma3.2} (1) follows from \eqref{eqn:(4.3)}.
This proves (1).
\par{\bf (2) }
By the definition of $g_{i}(t)$, we get
\begin{equation}\label{eqn:(3.4)}
\left(\frac{\partial_{t\bar{t}}\varphi_{i}}{\varphi_{i}}+\frac{\partial_{t\bar{t}}\varphi_{j}}{\varphi_{j}}\right)
\frac{g_{i}g_{j}}{g^{2}}
=
\frac{(\varphi_{i}\partial_{t\bar{t}}\varphi_{j}+\varphi_{j}\partial_{t\bar{t}}\varphi_{i})\cdot
|t|^{2(r_{i}+r_{j})}(\log|t|^{2})^{\ell_{i}+\ell_{j}}}{g^{2}},
\end{equation}
\begin{equation}\label{eqn:(3.5)}
\left(\frac{\partial_{t}\varphi_{i}}{\varphi_{i}}-\frac{\partial_{t}\varphi_{j}}{\varphi_{j}}\right)
\frac{g_{i}g_{j}}{g^{2}}
=
\frac{(\varphi_{j}\partial_{t}\varphi_{i}-\varphi_{i}\partial_{t}\varphi_{j})\cdot
|t|^{2(r_{i}+r_{j})}(\log|t|^{2})^{\ell_{i}+\ell_{j}}}{g^{2}},
\end{equation}
\begin{equation}\label{eqn:(3.6)}
\left(\frac{\partial_{t}\varphi_{i}}{\varphi_{i}}\cdot\overline{\frac{\partial_{t}\varphi_{j}}{\varphi_{j}}}\right)
\frac{g_{i}g_{j}}{g^{2}}
=
\frac{(\partial_{t}\varphi_{i}\cdot\overline{\partial_{t}\varphi_{j}})\cdot
|t|^{2(r_{i}+r_{j})}(\log|t|^{2})^{\ell_{i}+\ell_{j}}}{g^{2}}.
\end{equation}
Since the functions
$\varphi_{i}\partial_{t\bar{t}}\varphi_{j}+\varphi_{j}\partial_{t\bar{t}}\varphi_{i}$,
$\varphi_{j}\partial_{t}\varphi_{i}-\varphi_{i}\partial_{t}\varphi_{j}$,
$\partial_{t}\varphi_{i}\cdot\overline{\partial_{t}\varphi_{j}}$ are bounded near $t=0$ and since
$|t|^{2(r_{i}+r_{j})}(-\log|t|^{2})^{\ell_{i}+\ell_{j}}\leq(-\log|t|^{2})^{2N}$ by the definition of $N$, 
we get (2) by the second equality of Lemma~\ref{Lemma3.2} (1) and
\eqref{eqn:(3.4)}, \eqref{eqn:(3.5)}, \eqref{eqn:(3.6)}.
\end{pf}

\begin{lemma}\label{Lemma3.3}
Let $\varphi_{i}\in C^{\infty}_{\bf R}(T)$ for $0\leq i\leq N$ and set 
$g(t)=\sum_{i=0}^{N}(\log|t|^{2})^{i}\varphi_{i}(t)$.
Assume that $g(t)>0$ on $T^{o}$ and that $\varphi_{i}(0)\not=0$ for some $0\leq i\leq N$. Set 
$$
\ell:=\max_{0\leq i\leq N,\,\varphi_{i}(0)\not=0}\{i\}\in{\bf Z}_{\geq0}.
$$
Then there exists a constant $C>0$ such that the following inequalities hold
$$
|\partial_{t}\log g(t)|\leq\frac{C}{|t|(-\log|t|)},
\qquad
\left|
\partial_{t\bar{t}}\log g(t)
+
\frac{\ell}{|t|^{2}(-\log|t|)^{2}}
\right|
\leq
\frac{C}{|t|^{2}(-\log|t|)^{3}}.
$$
\end{lemma}

\begin{pf}
Set $I=\{0,1,\ldots,N\}$ and $g_{i}(t):=(-\log|t|)^{i}\varphi_{i}(t)$ for $i\in I$. 
Namely, we set $(r_{i},\ell_{i})=(0,i)$ in Lemma~\ref{Lemma3.2}.
Since $g(t)=\varphi_{\ell}(0)(-\log|t|)^{\ell}(1+O(1/\log|t|))$ as $t\to0$, 
we get for each $0\leq i\leq N$ the following asymptotic behavior as $t\to0$:
\begin{equation}\label{eqn:(3.7)}
\left|\frac{g_{i}(t)}{g(t)}\right|
=
\begin{cases}
\begin{array}{ll}
O(|t|(-\log|t|)^{i})&(i>\ell),
\\
1+O(|t|(\log|t|)^{n})&(i=\ell),
\\
O((-\log|t|)^{-(\ell-i)})&(i<\ell).
\end{array}
\end{cases}
\end{equation}
By the first equality of Lemma~\ref{Lemma3.2} (1) and \eqref{eqn:(3.7)}, there are constants $C,C'>0$ such that
$$
|\partial_{t}\log g(t)|
\leq
\frac{C}{|t|(-\log|t|)}\sum_{i=0}^{N}\left|\frac{g_{i}}{g}\right|
\leq
\frac{C'}{|t|(-\log|t|)}.
$$
This proves the first inequality. 
Since $g(t)=\varphi_{\ell}(0)(-\log|t|)^{\ell}(1+O(1/\log|t|))$, there exists $c>0$
such that $g(t)\geq c>0$ on $T^{o}$. In particular $O(1/g(t))=O(1)$.
This, together with Lemma~\ref{Lemma3.2} (2), yields that
\begin{equation}\label{eqn:(3.8)}
\begin{aligned}
\partial_{t\bar{t}}\log g
&=
-\frac{\ell}{|t|^{2}(\log|t|^{2})^{2}}
-
\frac{1}{2}\sum_{(i,j)\not=(\ell,\ell)}
\frac{i+j}{|t|^{2}(\log|t|^{2})^{2}}\left(\frac{g_{i}}{g}\right)\left(\frac{g_{j}}{g}\right)
\\
&\quad
+\frac{1}{2}\sum_{i\not=j}
\frac{(i-j)^{2}}{|t|^{2}(\log|t|^{2})^{2}}
\left(\frac{g_{i}}{g}\right)\left(\frac{g_{j}}{g}\right)
+
O\left(\frac{(-\log|t|)^{2N}}{|t|}\right).
\end{aligned}
\end{equation}
Since $|g_{i}(t)g_{j}(t)/g(t)^{2}|=O(1/\log|t|)$ when $i\not=j$ by \eqref{eqn:(3.7)}, 
the second and the third term in the right hand side of \eqref{eqn:(3.8)} is bounded by $|t|^{-2}(-\log|t|)^{-3}$
as $t\to0$. Similarly, it follows from \eqref{eqn:(3.7)} that
the second term of the right hand side of \eqref{eqn:(3.8)} is bounded by
$|t|^{-2}(-\log|t|)^{-3}$. The second inequality follows from \eqref{eqn:(3.8)}.
\end{pf}

\section
{ The boundary behavior of the curvature of the $L^{2}$-metric}\label{sect:4}
\par
In this section, we define $N,\ell_{q}\in{\bf Z}_{\geq0}$ as
$$
N:=n\rho_{q},
\qquad
\ell_{q}:=\max_{0\leq i\leq N,\,a_{i}(0)\not=0}\{i\},
$$
where $a_{i}(t)\in C^{\infty}(T)$, $0\leq i\leq N$, are the functions in Theorem~\ref{Theorem2.1}.
Recall that the integer $\nu>0$ was defined in Sect.\,\ref{sect:2}.

\subsection
{The singularity of the first Chern form}\label{subsect:4.1}
\par

\begin{theorem}\label{Theorem4.1}
The following formula holds as $s\to0$:
$$
c_{1}(R^{q}\pi_{*}\omega_{X/S}(\xi),h_{L^{2}})
=
\left\{
\frac{\ell_{q}}{|s|^{2}(\log|s|)^{2}}+O\left(\frac{1}{|s|^{2}(\log|s|)^{3}}\right)
\right\}
\sqrt{-1}\,ds\wedge d\bar{s}
$$
\end{theorem}

\begin{pf}
Recall that $T$ is another unit disc and that the map $\mu\colon T\to S$  is defined as $s=\mu(t)=t^{\nu}$. 
By Theorem~\ref{Theorem2.1}, we get
$$
\begin{aligned}
\mu^{*}c_{1}(R^{q}\pi_{*}\omega_{X/S}(\xi),h_{L^{2}})
&=
-\frac{\sqrt{-1}}{2\pi}\mu^{*}\partial\bar{\partial}\log\det G(s)
=
-\frac{\sqrt{-1}}{2\pi}\partial\bar{\partial}\log\det H(t)
\\
&=
-\frac{\sqrt{-1}}{2\pi}\partial\bar{\partial}\log[\sum_{m=0}^{N}a_{m}(t)\,(\log|t|^{2})^{m}].
\end{aligned}
$$
We set $g(t)=\det H(t)=\sum_{i=0}^{N}(\log|t|^{2})^{i}a_{i}(t)$ in Lemma~\ref{Lemma3.3}.
Since $a_{i}(0)\not=0$ for some $0\leq i\leq N$ by Theorem~\ref{Theorem2.1},
we deduce from Lemma~\ref{Lemma3.3} that
\begin{equation}\label{eqn:(4.6)}
\mu^{*}c_{1}(R^{q}\pi_{*}\omega_{X/S}(\xi),h_{L^{2}})
=
\ell_{q}\frac{\sqrt{-1}\,dt\wedge d\bar{t}}{|t|^{2}(-\log|t|)^{2}}
+
O\left(\frac{\sqrt{-1}\,dt\wedge d\bar{t}}{|t|^{2}(-\log|t|)^{3}}\right).
\end{equation}
Since 
$\mu^{*}\{\sqrt{-1}ds\wedge d\bar{s}/(|s|^{2}(-\log|s|)^{m}\}=
\nu^{2-m}\sqrt{-1}dt\wedge d\bar{t}/(|t|^{2}(-\log|t|)^{m})$, 
the desired inequality follows from \eqref{eqn:(4.6)}.
\end{pf}

\begin{remark}\label{Remark4.2}
The Hermitian metric $\mu^{*}\det h_{L^{2}}$ on the line bundle $\det R^{q}f_{*}\omega_{Y/T}(\xi)$ 
is good in the sense of Mumford \cite{Mumford77}. Namely, the following estimates hold: 
\begin{itemize}
\item[(1)]
There exist constants $C,\ell>0$ such that
$$
\det H(t)\leq C(-\log|t|^{2})^{\ell},
\qquad
(\det H(t))^{-1}\leq C(-\log|t|^{2})^{\ell}.
$$
\item[(2)]
There exists a constant $C>0$ such that
$$
|\partial_{t}\log\det H(t)|\leq\frac{C}{|t|(-\log|t|)},
\qquad
|\partial_{t\bar{t}}\log\det H(t)|\leq\frac{C}{|t|^{2}(-\log|t|)^{2}}.
$$
\end{itemize}
The inequalities (1) follow from Theorem~\ref{Theorem2.1}.
By setting $g(t)=\det H(t)$ in Lemma~\ref{Lemma3.3}, we get (2) because
$\det H(t)=g(t)=\sum_{i=0}^{N}(\log|t|^{2})^{i}a_{i}(t)$, $a_{i}(t)\in C^{\infty}_{\bf R}(T)$
with $a_{i}(0)\not=0$ for some $0\leq i\leq N$ by Theorem~\ref{Theorem2.1}.
\par
We do not know if the $L^{2}$-metric $\mu^{*}h_{L^{2}}$ on $R^{q}f_{*}\omega_{Y/T}(F^{*}\xi)$ is good 
in the sense of Mumford, because the estimates
$$
\|\partial_{t}H\cdot H^{-1}\|\leq\frac{C}{|t|(-\log|t|)},
\qquad
\|\partial_{\bar{t}}(\partial_{t}H\cdot H^{-1})\|\leq\frac{C}{|t|^{2}(-\log|t|)^{2}}
$$
do not necessarily follow from Theorem~\ref{Theorem2.1}; from Theorem~\ref{Theorem2.1}, 
we have only the estimates
$\|\partial_{t}H\cdot H^{-1}\|\leq C(-\log|t|)^{\ell}/|t|$ and
$\|\partial_{\bar{t}}(\partial_{t}H\cdot H^{-1})\|\leq C(-\log|t|)^{\ell}/|t|^{2}$,
where $\|A\|=\sum_{i,j}|a_{ij}|$ for a matrix $A=(a_{ij})$.
\end{remark}

\subsection
{Proof of Theorem~\ref{MainTheorem}}\label{subsect:4.2}
\par
Let $\lambda_{1},\ldots,\lambda_{\rho_{q}}$ be the eigenvalues of the Hermitian endomorphism 
$\sqrt{-1}{\mathcal R}(s)$.
By the Nakano semi-positivity of $(R^{q}\pi_{*}\omega_{X/S}(\xi),h_{L^{2}})$, we get $\lambda_{\alpha}\geq0$
for all $1\leq\alpha\leq\rho_{q}$. By Theorem~\ref{Theorem4.1}, we have the following inequality on $S^{o}$
$$
0\leq \sqrt{-1}\,{\rm Tr}[{\mathcal R}(s)]=\sum_{\alpha}\lambda_{\alpha}\leq\frac{C}{|s|^{2}(-\log|s|)^{2}}.
$$
In particular, we get $\Lambda:=\max_{\alpha}\{\lambda_{\alpha}\}\leq C/(|s|^{2}(-\log|s|)^{2})$.
We get the desired inequality for $\sqrt{-1}{\mathcal R}(s)$ from the inequality
$\sqrt{-1}{\mathcal R}(s)\leq\Lambda\cdot{\rm Id}_{R^{q}\pi_{*}\omega_{X/S}(\xi)}$. 
The inequality for $c_{1}(R^{q}\pi_{*}\omega_{X/S}(\xi),h_{L^{2}})$ is already proved
in Theorem~\ref{Theorem4.1}. This completes the proof.
\qed

\subsection
{Proof of Theorem~\ref{Theorem:HessianTorsion}}\label{subsect:4.3}
\par
By the curvature formula for Quillen metrics \cite{BGS88}, 
the following equation of currents on $S^{o}$ holds
\begin{equation}\label{eqn:(4.7)}
\begin{aligned}
\,&
-dd^{c}\log\tau(X/S,\omega_{X/S})+\sum_{q\geq0}(-1)^{q}c_{1}(R^{q}\pi_{*}\omega_{X/S}(\xi),h_{L^{2}})
\\
&=
[\pi_{*}\{{\rm Td}(TX/S,h_{X/S}){\rm ch}(\omega_{X/S}(\xi))\}]^{(2)},
\end{aligned}
\end{equation}
where $[A]^{(p)}$ denotes the component of degree $p$ of a differential form $A$.
By \cite[Lemma 9.2]{Yoshikawa07}, there exists $r\in{\bf Q}_{>0}$ such that as $s\to0$
\begin{equation}\label{eqn:(4.8)}
[\pi_{*}\{{\rm Td}(TX/S,h_{X/S}){\rm ch}(\omega_{X/S}(\xi))\}]^{(2)}(s)
=
O\left(\frac{\sqrt{-1}\,|s|^{2r}(-\log|s|)^{n}ds\wedge d\bar{s}}{|s|^{2}}\right).
\end{equation}
By Theorem~\ref{MainTheorem}, we get
\begin{equation}\label{eqn:(4.9)}
\sum_{q\geq0}(-1)^{q}c_{1}(R^{q}\pi_{*}\omega_{X/S}(\xi),h_{L^{2}})
=
\frac{\sum_{q\geq0}(-1)^{q}\ell_{q}}{2\pi}\frac{\sqrt{-1}\,ds\wedge d\bar{s}}{|s|^{2}(-\log|s|)^{2}}
+
O\left(\frac{\sqrt{-1}\,ds\wedge d\bar{s}}{|s|^{2}(-\log|s|)^{3}}\right).
\end{equation}
By \eqref{eqn:(4.7)}, \eqref{eqn:(4.8)}, \eqref{eqn:(4.9)}, we get on $S^{o}$
$$
\frac{\sqrt{-1}}{2\pi}\partial\bar{\partial}\log\tau(X/S,\omega_{X/S})
=
\frac{\sum_{q\geq0}(-1)^{q}\ell_{q}}{2\pi}\,\frac{\sqrt{-1}\,ds\wedge d\bar{s}}{|s|^{2}(-\log|s|)^{2}}
+
O\left(\frac{\sqrt{-1}\,ds\wedge d\bar{s}}{|s|^{2}(-\log|s|)^{3}}\right).
$$
This completes the proof.
\qed

\section
{Canonical singularities and the curvature of $L^{2}$-metric}\label{sect:5}
\par
In this section, we assume that the central fiber $X_{0}$ is reduced and irreducible and has only canonical 
(equivalently rational) singularities. Then $G(s)=(G_{\alpha\bar{\beta}}(s))$ is expected to have
better regularity than usual. To see this, set
$$
{\mathcal B}(S)
:=
C^{\infty}(S)\oplus\bigoplus_{r\in{\bf Q}\cap(0,1]}\bigoplus_{k=0}^{n}|s|^{2r}(\log|s|)^{k}C^{\infty}(S)
\subset
C^{0}(S).
$$
By \cite[Th.\,7.2]{Yoshikawa10}, the $L^{2}$-metric $h_{L^{2}}$ on $R^{q}\pi_{*}\omega_{X/S}(\xi)$ 
is a continuous Hermitian metric lying in the class ${\mathcal B}(S)$. 
Namely, $G_{\alpha\bar{\beta}}(s)\in{\mathcal B}(S)$ for all $1\leq\alpha,\beta\leq\rho_{q}$.

\begin{proposition}\label{Proposition5.1}
If $X_{0}$ has only canonical singularities, then
there exists $r\in{\bf Q}_{>0}$ and $C>0$ such that the following inequality of real $(1,1)$-forms on $S^{o}$
holds
$$
0
\leq 
c_{1}(R^{q}\pi_{*}\omega_{X/S}(\xi),h_{L^{2}})
\leq 
C\frac{\sqrt{-1}\,|s|^{2r}ds\wedge d\bar{s}}{|s|^{2}(-\log|s|)^{2}}.
$$
In particular, the curvature $i{\mathcal R}(s)\,ds\wedge d\bar{s}$ satisfies the following estimate:
$$
0\leq \sqrt{-1}{\mathcal R}(s)\leq\frac{C|s|^{2r}}{|s|^{2}(-\log|s|)^{2}}{\rm Id}_{R^{q}\pi_{*}\omega_{X/S}(\xi)}.
$$
\end{proposition}

\begin{pf}
Since $G_{\alpha\bar{\beta}}(s)$ is continuous on $S$, we may assume by an appropriate choice of basis
that $G_{\alpha\bar{\beta}}(0)=\delta_{\alpha\beta}$.
Since $\det G(s)\in{\mathcal B}(S)$, there exist a finite set $I$ and
$(r_{i},\ell_{i})\in{\bf Q}_{>0}\times{\bf Z}_{\geq0}$ for each $i\in I$ such that
$$
\det G(s)=1+\sum_{i\in I}|s|^{2r_{i}}(\log|s|^{2})^{\ell_{i}}\,\varphi_{i}(s).
$$
We set $r_{0}=0$, $\ell_{0}=0$ and $\varphi_{0}(s)=1$. For $i\in\{0\}\cup I$, we set
$g_{i}(s):=|s|^{2r_{i}}(\log|s|^{2})^{\ell_{i}}\,\varphi_{i}(s)$.
By Lemma~\ref{Lemma3.2} (2) applied to $\det G(s)$, we get
$$
\begin{aligned}
\,&
-\partial_{s\bar{s}}\log\det G(s)
\\
&=
\frac{1}{2}\sum_{i,j\in I\cup\{0\}}\frac{\ell_{i}+\ell_{j}}{|s|^{2}(\log|s|^{2})^{2}}\cdot\frac{g_{i}g_{j}}{g^{2}}
-
\frac{1}{2}\sum_{i,j\in I\cup\{0\}}
\left|
\frac{r_{i}-r_{j}}{s}+\frac{\ell_{i}-\ell_{j}}{s(\log|s|^{2})}
\right|^{2}
\frac{g_{i}g_{j}}{g^{2}}
\\
&\qquad
+O\left(\frac{(-\log|t|)^{2N}}{|t|}\right)
\\
&\leq
C\sum_{i,j\in I\cup\{0\}}\frac{(\ell_{i}+\ell_{j})|s|^{2(r_{i}+r_{j})}}{|s|^{2}(\log|s|^{2})^{2}}
+
C\sum_{i,j\in I\cup\{0\}}\left|\frac{r_{i}-r_{j}}{s}+\frac{\ell_{i}-\ell_{j}}{s(\log|s|^{2})}\right|^{2}|s|^{2(r_{i}+r_{j})}
\\
&\qquad
+O\left(\frac{(-\log|t|)^{2N}}{|t|}\right).
\end{aligned}
$$
Set $r:=\min_{i\in I}\{r_{i}\}>0$.
Since $r_{i}+r_{j}>r$ for all $(i,j)\in(I\cup\{0\})\times(I\cup\{0\})\setminus\{(0,0)\}$, we get
$$
-\partial_{s\bar{s}}\log\det G(s)
\leq
C\frac{2\ell_{0}}{|s|^{2}(\log|s|^{2})^{2}}+C\sum_{(i,j)\not=(0,0)}\frac{|s|^{2(r_{i}+r_{j})}}{|s|^{2}(\log|s|^{2})^{2}}
\leq
\frac{C\,|s|^{2r}}{|s|^{2}(\log|s|^{2})^{2}}.
$$
because $\ell_{0}=0$. Since $-\partial_{s\bar{s}}\log\det G(s)\geq0$ by the Nakano semi-positivity of
$(R^{q}\pi_{*}\omega_{X/S}(\xi),h_{L^{2}})$ by \cite{Berndtsson09}, \cite{MourouganeTakayama08},
we get the first inequality.
The proof of the second inequality is the same as that of the corresponding inequality of 
Theorem~\ref{MainTheorem}.
\end{pf}


\end{document}